\documentclass[11pt,a4paper]{article}
\usepackage{amsmath}
\usepackage{amssymb}
\usepackage{amsthm}
\usepackage{latexsym}
\usepackage{amsbsy}
\usepackage[all]{xypic}
\usepackage{amscd}

\newtheorem{thm}{Theorem}[section]
\newtheorem{cor}[thm]{Corollary}
\newtheorem{lem}[thm]{Lemma}

\newtheorem{prop}[thm]{Proposition}
\newtheorem{defn}[thm]{Definition}

\begin{document}

\begin{center}
{\Large \bf
From triangulated categories to abelian categories \\ --  cluster
tilting in a general framework}
\bigskip

{\large Steffen Koenig and
 Bin Zhu\footnote{Supported by the NSF of China (Grants 10471071) and
by The Leverhulme Trust through the network 'Algebras, Representations
and Applications'.}}
\bigskip

{\small
\begin{tabular}{cc}
Mathematisches Institut & Department of Mathematical Sciences
\\
Universit\"at zu K\"oln & Tsinghua University
\\
Weyertal 86-90, 50931 K\"oln, Germany &   100084 Beijing, P. R. China
\\
{\footnotesize E-mail: skoenig@math.uni-koeln.de} &
{\footnotesize E-mail: bzhu@math.tsinghua.edu.cn}
\end{tabular}
}
\bigskip


\end{center}

\def\s{\stackrel}
\def\gama{\gamma}
\def\Longrightarrow{{\longrightarrow}}
\def\P{{\cal P}}
\def\A{{\cal A}}
\def\F{\mathcal{F}}
\def\X{\mathcal{X}}
\def\T{\mathcal{T}}
\def\m{\textbf{ M}}
\def\t{{\tau }}
\def\b{\textbf{d}}
\def\K{{\cal K}}

\def\G{{\Gamma}}
\def\e{\mbox{exp}}

\def\righta{\rightarrow}

\def\s{\stackrel}

\def\ncong{\not\cong}

\def\mathbb{\NN}

\def\Hom{\mbox{Hom}}
\def\Ext{\mbox{Ext}}
\def\ind{\mbox{ind}}
\def\coprod{\amalg }
\def\L{\Lambda}
\def\c{\circ}
\def\mu{\multiput}

\renewcommand{\mod}{\operatorname{mod}\nolimits}
\newcommand{\add}{\operatorname{add}\nolimits}
\newcommand{\Rad}{\operatorname{Rad}\nolimits}
\newcommand{\RHom}{\operatorname{RHom}\nolimits}
\newcommand{\uHom}{\operatorname{\underline{Hom}}\nolimits}
\newcommand{\End}{\operatorname{End}\nolimits}
\renewcommand{\Im}{\operatorname{Im}\nolimits}
\newcommand{\Ker}{\operatorname{Ker}\nolimits}
\newcommand{\Coker}{\operatorname{Coker}\nolimits}
\renewcommand{\r}{\operatorname{\underline{r}}\nolimits}
\def \text{\mbox}

\hyphenation{ap-pro-xi-ma-tion}

\begin{abstract}
A general framework for cluster tilting is set up by showing that any quotient
of a triangulated category modulo a tilting subcategory (that is, a maximal
one-orthogonal subcategory) carries an induced abelian structure. These
abelian quotients turn out to be module categories of Gorenstein algebras
of dimension at most one.

\end{abstract}

\begin{center}


\medskip

\end{center}
\medskip

\section{Introduction}

Abelian and triangulated categories are two fundamental structures
in algebra and geometry. While modules or sheaves are forming
abelian categories, complexes lead to homotopy or derived categories
that are triangulated. Triangulated categories which at the same
time are abelian must be semisimple. There are, however, well-known
ways to produce triangulated categories from abelian ones [15]. For
example, taking the module category of a self-injective algebra
'modulo projectives' produces the stable category, which is
triangulated. In particular, \cite{H} the (triangulated) derived
module category of a finite dimensional algebra of finite global
dimension is equivalent to the stable category of its (infinite
dimensional, but locally finite-dimensional)  repetitive algebra.
Homotopy categories of complexes provide another example of passing
from abelian categories (of complexes) to triangulated ones.

Among the surprises produced by the recent theory of cluster algebras and
cluster categories is the possibility of sometimes going the opposite way;
starting from a cluster category, which is a triangulated category
constructed from a derived category, one can pass to a quotient category,
which turns out to be abelian \cite{BMR,BMRRT,Z1,KR}.
The quotient is taken modulo a 'tilting subcategory' (a maximal 1-orthogonal
subcategory as defined in \cite{I1, I2}). Such tilting subcategories
correspond to clusters in the cluster algebras introduced by Fomin and
Zelevinsky \cite{FZ1,FZ2}. Their endomorphism rings have been shown to have
interesting properties \cite{KR}, such as being
Gorenstein of dimension one, and are expected to contain crucial
information about clusters and cluster variables.

In this article we are going to provide a general framework for passing
from triangulated categories to abelian categories by factoring out
tilting subcategories. Indeed, our main result states in full generality
that any such quotient category carries an abelian structure.
We will relate the two structures in a direct and
explicit way, thus not only reproving, but also strengthening the known
results in the case of cluster categories. In particular, we give
explicit constructions of kernels and cokernels in the abelian quotient
category. By examples we show that our result also applies to stable
categories and that it is not restricted to Calabi Yau dimension two.

A crucial difference between abelian and triangulated categories concerns
monomorphisms and epimorphisms. In an abelian category, plenty of these
must exist, as kernels or cokernels, respectively. In a triangulated
category, however, any monomorphism or epimorphism (in the categorical
sense) must split. One of our main tools is a characterisation of maps
in a triangulated category which become monomorphisms or epimorphisms after
factoring out a tilting subcategory.

The abelian quotient categories in our general situation enjoy
nice properties similar to the situation of Calabi Yau dimension two studied
by Keller and Reiten \cite{KR}. For instance, they have enough projectives
and injectives and thus are equivalent to module categories over the
endomorphism rings of tilting objects. Moreover, these endomorphism
rings are Gorenstein of dimension at most one.
\bigskip

The article is organised as follows:

In Section 2 we first collect basic
material on quotient categories and then prove Theorem \ref{charmonoepi},
the characterisation of morphisms in a triangulated category which
become monomorphisms or epimorphisms in a quotient modulo a tilting
subcategory.

Section 3 contains the main result, Theorem \ref{mainthm},
and its proof; the quotient of a triangulated category
modulo a tilting subcategory carries an induced abelian structure.

In Section 4 we first extend a number of results on
the abelian quotient category, which are known for cluster categories or
Calabi-Yau categories of CY-dimension two, to quotients of triangulated
categories. In particular, we show in Theorem \ref{Gorenstein} that the
abelian quotient category always is the module category of a Gorenstein
algebra (maybe of infinite dimension or just a ring) of
Gorenstein dimension at most one. Moreover, we give examples different from
cluster categories or Calabi-Yau categories. Then we go back from our
general  cluster-tilting to the classical cluster-tilting theory and show
that the tilting subcategories of $D^b(H)$ are sent to cluster-tilting
subcategories (or cluster-tilting objects) by the projection
$\pi :D^b(H)\longrightarrow \mathcal{C}(H)$, thus complementing
results in \cite{BMR, BMRRT}. Moreover we show that the projection $\pi$
gives a covering functor from the subcategory of projective objects in the
abelian quotient $D^b(H)/\mathcal{T}$ to the
subcategory of projective modules of  the module category over the
corresponding cluster-tilting subcategory (cluster-tilted algebra)
$\pi(\mathcal{T})$ and that it also gives the corresponding push-down
functor between their module categories; this again accompanies results in
\cite{BMRRT}.

Section 5 discusses various aspects of potential converses of Theorem
\ref{mainthm}. There are trivial and non-trivial counterexamples to a
direct converse. Assuming an abelian structure on a quotient of a
triangulated category modulo some subcategory, we can, however, recover
some of the conditions used as assumptions in Theorem \ref{mainthm} and
also the characterisation in Theorem \ref{charmonoepi}.
\bigskip

\section{Quotient categories and morphisms}

In this section, we collect some basic material and then prove our main
tool, Theorem \ref{charmonoepi}, characterising morphisms which become
monomorphisms or epimorphisms in a quotient category.

\subsection{Basics on quotient categories}

Let $\cal H$ be an additive category and $\cal T$ a full subcategory which
is closed under taking direct sums and direct summands, i.e, for any two
objects $X, Y\in \mathcal{H},\  X\oplus Y\in \mathcal{T}$ if and only
if $X,\ Y \in \mathcal{T}$.
Then the
quotient category ${\cal A}:={\cal H}/{\cal T}$ has the same objects as
$\cal H$, and morphisms from $X$ to $Y$ in the quotient category are
the $\cal H$-morphisms modulo the subgroup of morphisms factoring through
some object in $\mathcal T$. For $f$ a morphism in $\mathcal H$,
we denote by $\underline{f}$ its residue class in the quotient category.
The quotient $\mathcal{A}$ is
additive. If $\mathcal{T }=\mbox{add}T$ for some object $T$, the quotient
category is denoted by $\mathcal{H}/T $.

The following is well-known:

\begin{lem}
(a) The property of being a Krull-Schmidt category is inherited by the
quotient category.

(b) $\ind\mathcal{A}=ind\mathcal{H}\setminus ind\mathcal{T}$.
\end{lem}

Throughout the paper, the shift functor of a triangulated category will be
denoted by $[1]$.
 From Subsections 4.2 to the end of the paper,  we will assume that
$\mathcal{H}$
is a $k-$linear triangulated category with split idempotents. Furthermore,
in these sections we also
will assume that all Hom-spaces of $\mathcal{H}$ are finite
dimensional and the existence of a Serre functor $\Sigma$ on $\mathcal H$,
that is, an autoequivalence naturally satisfying
$Hom(X, Y) \simeq DHom(Y,\Sigma X)$.
Hence, $\mathcal H$ has Auslander-Reiten triangles, and
$\Sigma=\tau [1]$. Here,
$\tau$ is the Auslander-Reiten translate. The space
$Hom(Y,X[1])$ sometimes is denoted by $Ext^1(Y,X)$. If $Ext^1(T,T')=0$
for any $T, T'\in \mathcal{T}$, we
say that $\mathcal{T}$ satisfies $Ext^1(\mathcal{T},\mathcal{T})=0.$

If $\cal H$ is a triangulated category, then its distinguished triangles
will be just called triangles.

\subsection{How to become a monomorphism}

A morphism $f$ in a category $\mathcal{H}$ is called
a {\em monomorphism} (or an
{\em epimorphism}) provided $g=0$ whenever $f \circ g=0$ (respectively
$g \circ f=0$).

The following well-known lemma exhibits a crucial difference between
triangulated and abelian categories.

\begin{lem}
A monomorphism in a triangulated category is a section, that is, it
admits a left inverse. Dually an epimorphism in a triangulated
category is a retraction, that is, it admits a right inverse.
\end{lem}

\begin{proof}
See, for example, Exercise 1 in \cite{GM}, IV.1.
\end{proof}

\begin{thm} \label{charmonoepi}
Let $\mathcal{H}$ be a triangulated category and $\mathcal{T }$ a full
subcategory of
$\mathcal{H}$ with $Ext^1(\mathcal{T},\mathcal{T})=0$.
Let $f:X\rightarrow Y$ be a morphism in $\mathcal{H}$ which is
a part of a triangle
$Z[-1]\s{h}{\rightarrow}X\s{f}{\rightarrow}Y\s{g}{\rightarrow}Z$.
Then $\underline{f}$ is a monomorphism in $\mathcal{A}$ if and only
$\underline{h}=0;$  $\underline{f}$ is an epimorphism in $\mathcal{A}$
if and only if $\underline{g}=0.$

In particular, if $Z$ is in  $\mathcal{T }$, then $\underline{f}$ is
an epimorphism; if $Z[-1]$ is in  $\mathcal{T }$, then $\underline{f}$ is
a monomorphism.
\end{thm}

\begin{proof}
We give the proof for the statement about epimorphisms, the case of
monomorphisms being dual.
We first will deal with a special case:
Suppose that $f:X\rightarrow Y$ is a morphism in $\mathcal{H}$ which is
part of a triangle
$T[-1]\s{h}{\rightarrow}X\s{f}{\rightarrow}Y\s{g}{\rightarrow}T$
with $T\in \mathcal{T }.$
We will prove that $\underline{f}$ is an
epimorphism in $\mathcal{A}$.
Let $g_1: Y\rightarrow Z$ with
$\underline{g_1} \circ \underline{f}=0.$ Then there exists an object
$T'\in\mathcal{T }$ and morphisms $g'_1:X\rightarrow T'$ and
$f'_1:T'\rightarrow Z$ with $g_1 \circ f=f'_1 \circ g'_1$.
The resulting commutative
square can be completed to a commutative diagram with rows being triangles:

\[ \begin{CD}
T[-1]@>h>>X @>f>>Y@>g>>T\\
@VVV@VVg'_1 V @VVg_1V @VV V \\
M[-1]@>>>T'@>f'_1>>Z@>>>M
    \end{CD} \]

In this diagram, $g'_1 \circ h=0$ since
${\Hom}(T[-1], T')\cong\mbox{Hom}(T,T'[1])=0$.
By the long exact homology sequence associated to the triangle,
the morphism $g'_1$ factors through $Y$, i.e. there
is a morphism $f_1:Y\rightarrow T'$ with $g'_1=f_1 \circ f$. It follows
that $(g_1-f'_1 \circ f_1) \circ f=0$, and then the morphism
$g_1-f'_1 \circ f_1$ factors through
$g$, i.e., there is a morphism $\sigma_1: T\rightarrow Z$ such that
$g_1-f'_1 \circ f_1=\sigma _1 \circ g.$ Therefore we have
$g_1=f'_1\circ f_1+\sigma _1 \circ g$
which means that $g_1$ factors through $T\oplus T'$. Hence
$\underline{g_1}=0,$ proving that  $\underline{f}$ is epimorphism.

Now we turn to the general case: Suppose that $f:X\rightarrow Y$
is a morphism in $\mathcal{H}$ which is part of a triangle
$Z[-1]\s{h}{\rightarrow}X\s{f}{\rightarrow}Y\s{g}{\rightarrow}Z$
such that $g$ factors through $T$ with $T\in \mathcal{T}$. Then,
by completing the right hand square, we get
a commutative diagram with rows being triangles:

\[ \begin{CD}
T[-1]@>>>M@>f_1>>Y@>\sigma_1>>T\\
@VVV@VVg_1 V @VVidV @VV\sigma'_1 V \\
X[-1]@>h>>X @>f>>Y@>g>>Z   \end{CD} \] As shown above, we have that
$\underline{f_1}$ is an epimorphism in $\mathcal{A}$. This implies
that $\underline{f}$ is also an epimorphism in $\mathcal{A}$ since
$\underline{f} \circ \underline{g_1}=\underline{f_1}.$

The converse is easy: Suppose $\underline{f}$ is an epimorphism and
$\underline{g} \circ \underline{f} = 0$, then $\underline{g} =0$.
\end{proof}

\medskip

\begin{cor}
Let $f:X\rightarrow Y$ be a morphism in
$\mathcal{H}$ which is part of a triangle
$Z[-1]\s{h}{\rightarrow}X\s{f}{\rightarrow}Y\s{g}{\rightarrow}Z$. If
$\underline{f}$ and $\underline{g}$ are zero maps, then $Y\in
\mathcal{T}.$
\end{cor}

\begin{proof}
As $\underline{g}=0,$ $\underline{f}$ is an
epimorphism. Then $\underline{id}_Y \circ \underline{f} = 0$ implies
$\underline{id}_Y = 0$ and thus $Y\cong 0$ in $\mathcal{A}$, i.e.
$Y\in\mathcal{T}.$
\end{proof}

\section{Induced abelian structure on quotient categories}

In this section we prove our main result, Theorem \ref{mainthm}, stating
that any quotient of any triangulated category modulo any tilting
subcategory carries an induced abelian structure.

The following definition is due to Iyama \cite{I1}.

\begin{defn} \label{defmaxorth}
Let $\mathcal{H}$ be an abelian category or a triangulated category.
A subcategory $\mathcal{T }$ of $\mathcal{H}$ is called a {\em maximal
$1$-orthogonal subcategory} of $\mathcal{H}$,
if it satisfies the following conditions:

\begin{enumerate}

\item $\mathcal{T }$ is contravariantly finite and covariantly finite,

\item $X\in \mathcal{T }$ if and only if $\Ext^1(X,\mathcal{T })=0,$

\item $X\in \mathcal{T }$ if and only if $\Ext^1(\mathcal{T },X)=0.$

\end{enumerate}

In case $\mathcal{H}$ is a triangulated category, the maximal
$1$-orthogonal subcategories are called {\em tilting subcategories}.
If $\mathcal{T}=addT$, we call $T$ is {\em  a maximal $1-$orthogonal object}
of $\mathcal{H}.$
\end{defn}

Recall that a subcategory $\mathcal{T}$ is called contravariantly finite
in $\mathcal{H}$ provided for any object $X$ of $\mathcal{H}$ there is
a right $\mathcal{T}$-approximation $f: T\rightarrow X$, i.e. $Hom(T',f):
Hom(T',T)\longrightarrow Hom(T',X)$ is surjective for any $T'\in
\mathcal{T}.$ Dually, one can define left $\mathcal{T}-$approximation of
$X$ and covariantly finiteness of $\mathcal{T}$.

We first show that in a triangulated category, some of the defining
conditions of maximal 1-orthogonal imply others.

\begin{lem} \label{1-orthconditions}
Let $\mathcal{H}$ be a triangulated category and $\mathcal{T}$ a full
subcategory. Then
\begin{enumerate}
\item If $\mathcal{T}$ is contravariantly finite in $\mathcal{H}$ and
satisfies condition $3$ of Definition \ref{defmaxorth}, then for any
object $M$ of $\mathcal{H}$,
there is a triangle $T_1\rightarrow T_0\s{\sigma}{\rightarrow}
M\rightarrow T_1[1]$ with  $\sigma$ being a right
$\mathcal{T}-$approximation of $M$, and
 $T_1\in \mathcal{T}.$  Dually, if  $\mathcal{T}$ is covariantly finite
in $\mathcal{H}$ and satisfies condition $2$ of Definition \ref{defmaxorth},
 then  for any object $M$ of $\mathcal{H}$,
there is a triangle $M\s{\sigma}{\rightarrow} T_0\rightarrow
T_1\rightarrow M[1]$ with  $\sigma$ being a left
$\mathcal{T}-$approximation
of $M$, and  $T_1\in \mathcal{T}.$
\item Let $G$ be an automorphism of $\mathcal{H}$. Then $\mathcal{T}$ is
contravariantly (or covariantly) finite in $\mathcal{H}$ if and only if so is
the image $G(\mathcal{T})$.
\item If $\mathcal{T}$ is contravariantly finite and satisfies condition
$3$ of Definition \ref{defmaxorth}, then $\mathcal{T}$ is a
tilting subcategory of $\mathcal{H}$, i.e. it satisfies all conditions in
Definition \ref{defmaxorth}.
\end{enumerate}
\end{lem}

\begin{proof}
(1) We will prove the statement on right $\mathcal{T}-$approximations, the
case of a left $\mathcal{T}-$approximation being dual. Let $\sigma:\
T_0\rightarrow M$ be a right $\mathcal{T}-$approximation of $M$, and
$X\rightarrow T_0\s{\sigma}{\rightarrow} M\rightarrow X[1]$ a triangle
containing $\sigma$.
For any $T\in \mathcal{T}$, by applying Hom$(T,-)$ to
this triangle, there is an exact sequence $\mbox{Hom}(T,X)
\longrightarrow \mbox{Hom}(T, T_0)\s{Hom(T,\sigma )}{\longrightarrow}
 \mbox{Hom}(T, M)\longrightarrow \mbox{Hom}(T,X[1])\rightarrow
\mbox{Hom}(T, T_0[1])$. Since Hom$(T,\sigma )$ is surjective and
Hom$(T,T_0[1])=0$, also
Hom$(T,X[1])=0$. By condition $3$ it follows that $X\in \mathcal{T}$.

(2) We will prove that $G$ sends any right $\mathcal{T}-$approximation
$f:T\rightarrow M$ to a right $G(\mathcal{T})-$approximation of $G(M)$.
Given any morphism $g: G(T')\rightarrow G(M)$, we can write $g=G(g')$ where
$g': T'\rightarrow M$ since $G$ is full. Hence there is a map
$h:T'\rightarrow T$ with $g'=f \circ h$. Then $g=G(g')=G(f) \circ G(h)$.
Thus $G(f)$ is a right $G(\mathcal{T})-$approximation of $G(M)$. Since $G$ is
 an automorphism, any object of $\mathcal{H}$ is isomorphic to $G(X)$ for
some object $X$. Therefore the image of  $\mathcal{T}$ under $G$ is
contravariantly finite in $\mathcal{H}.$

(3) Suppose $\mathcal{T}$ is contravariantly finite in $\mathcal{H}$ and
satisfies condition 3 of Definition \ref{defmaxorth}. Then
$\mathcal{T}[-1]$ is also contravariantly finite
 by (2). Similarly the subcategory $\mathcal{T}[-1]$ satisfies an analogue
of condition $3$ of Definition \ref{defmaxorth} since $\mathcal{T}$
satisfies condition $3$. Then for any
$M\in \mathcal{H}$, it follows from assertion $1$ of this lemma
that there is a triangle   $T_1[-1]\rightarrow
T_0[-1]\s{\sigma}{\rightarrow} M\s{\beta}{\rightarrow} T_1$ with
$\sigma$ being right $\mathcal{T}[-1]-$approximation of $M$. Since
$Hom(T_0[-1],T)\cong Ext^1(T_0, T)=0$ for any $T\in \mathcal{T},$
 we have that $\beta$ is a left $\mathcal{T}-$approximation of $M$. This
proves that $\mathcal{T}$ is also covariantly finite in $\mathcal{H}$. Now
assume that $Ext^1(X,\mathcal{T})=0$ for some $X$.
We have to prove $X\in \mathcal{T}.$
Let $T_0[-1]\rightarrow X$ be the right $\mathcal{T}[-1]-$approximation of
$X$. Then we have a triangle $T_1 \rightarrow T_0 \rightarrow
X\s{h}{\rightarrow} T_1[1]$ by the statement of part (1).
Then $h\in Ext(X, \mathcal{T})=0.$  Thus
the triangle splits, $T_0\cong X\oplus T_1$ and thus $X\in
\mathcal{T}.$
\end{proof}

\begin{thm} Let $\mathcal{H}$ be a triangulated category
and $\mathcal{T }$ a tilting subcategory of $\mathcal{H}$. Then
$\mathcal{A} =\mathcal{H} / \mathcal{T } $ is an abelian category.
\label{mainthm}
\end{thm}

In the proof we will explicitly construct the abelian structure of
$\mathcal{A}$, that is, kernels and cokernels, from the triangulated
structure of $\mathcal{H}$.

\begin{proof} Since $\mathcal{A}$ is an additive category, in
order to prove it is abelian, we need to prove the existence of kernels
and cokernels and also that monomorphisms are kernels and epimorphisms are
cokernels.

{\bf Claim (1).}
For any morphism $f: X\rightarrow Y$, there is
a morphism $g: Y\rightarrow M$ which is the cokernel of $f$.

We complete $f$: $X\s{f}{\rightarrow} Y\s{f_1}{\rightarrow}
Z\s{f_2}{\rightarrow}X[1]$ to a triangle. Let $\sigma :
T_0\rightarrow X[1]$ be the right $\mathcal{T }-$approximation of
$X[1]$. Then we form another triangle $T_1\s{}{\rightarrow}
T_0\s{\sigma}{\rightarrow}X[1]\s{}{\rightarrow}T_1[1]$. Here
$\sigma$ being an approximation implies
that $T_1\in \mathcal{T }$. Composing the map $X[1] \rightarrow T_1[1]$
with $f_2$ we get a map $Z \rightarrow T_1[1]$. Extending to a triangle
we get the following commutative diagram:

$$ \begin{array}{cccclcl}
&&&&T_1&=&T_1\\
&&&&\downarrow&&\downarrow\\
&&Y&\s{g}{\rightarrow}&M&\rightarrow&T_0\\
&&\parallel&&\downarrow\sigma '&&\downarrow \sigma\\
X&\s{f}{\rightarrow}&Y&\s{f_1}{\rightarrow}&Z&\s{f_2}{\rightarrow}&X[1]
\end{array}$$

Then $\underline{f_1}=\underline{\sigma }' \circ \underline{g}$ and
$\underline{\sigma }'$ is a monomorphism, $\underline{g}$ is an
epimorphism by \ref{charmonoepi}.
Since $\sigma' \circ g \circ f = f_1 \circ f = 0$, we also have that
$\underline{g} \circ \underline{f}=0.$

We will prove that $\underline{g}$ is the cokernel of $\underline{f}.$

First, for any $h:Y\rightarrow N$ with $\underline{h} \circ
\underline{f}=0$ we will prove that $\underline{h}$ factors through
$\underline{g}$.

By $\underline{h} \circ \underline{f}=0$, it follows that $hf$ factors
through some object $T\in \mathcal{T }.$ Hence, there is the following
commutative diagram:

$$ \begin{array}{cclclcl}
&&Y&\s{g}{\rightarrow}&M&\rightarrow&T_0\\
&&\|&&\downarrow \sigma '&&\downarrow \sigma\\
X&\s{f}{\rightarrow} &Y&\s{f_1}{\rightarrow}&Z&\rightarrow&X[1]\\
\downarrow&&\downarrow h&&\downarrow h'&&\downarrow h''\\
T&\s{\rho}{\rightarrow}&N&\s{\rho _1}{\rightarrow}&Z'&\s{\rho
_2}{\rightarrow}&T[1]
\end{array}$$

Since $h''\circ \sigma\in \mbox{Hom}(T_0, T[1])=0$, we have that
$\rho_2 \circ h' \circ \sigma'=0$.
It follows that $h' \circ \sigma'$ factors through
$\rho_1$, i.e. there is a morphism $\rho':M\rightarrow N$ such that
$h' \circ \rho'=\rho_1 \circ \rho'$.
Therefore, $\rho_1 \circ (h-\rho'g)=0$ and $h-\rho' \circ g$ factors through
$\rho$, i.e. there is a morphism $\rho'': Y\rightarrow T$ such that
$h-\rho' \circ g=\rho \circ \rho''.$
So, $\underline{h}$ factors through $\underline{g}$, i.e.
$\underline{h}=\underline{\rho}' \circ \underline{g}$.

Second, we will prove that the map $\underline{\rho}'$ is unique.
Suppose that we have two such maps $\underline{\rho _1}'$ and
$\underline{\rho}'$ such that
$\underline{h}=\underline{\rho}' \circ \underline{g}=\underline{\rho
_1}' \circ \underline{g}$. Then $(\rho_1'-\rho') \circ g$ factors through some
object $T_1\in \mathcal{T }$, i.e.
$(\rho_1'-\rho') \circ g=\beta \circ \alpha$
with $\alpha: Y\rightarrow T_1$ and $\beta: T_1\rightarrow N$.
Let $T_1\s{\beta}{\rightarrow}
N\s{\beta_1}{\rightarrow} N'\s{\beta_2}{\rightarrow}T_1[1]$ be a
triangle into which $\beta$ is embedded. We have the commutative diagram:

\[ \begin{CD}
Y@>g>> M@>g'>>T@>>>Y[1]\\
@V\alpha VV@VV\rho_1'-\rho' V @VV\beta_3V @VV V \\
T_1@>\beta >>N@>\beta _1>>N'@>\beta _2>>T_1[1]
    \end{CD} \]

Since $\beta _2 \circ \beta_3\in \mbox{Hom}(T,T[1])=0$, $\beta _3$ factors
through $\beta_1$, i.e. $\beta _3=\beta_1 \circ \beta_4$ for a map
$\beta_4$. Therefore $\rho_1'-\rho'- \beta_4 \circ g'$ factors through
$\beta$. Then $\underline{\rho'}-\underline{\rho_1'}=0$. This
finishes the proof that any morphism in $\mathcal{A}$ has a unique cokernel.
Dually, we have that any morphism in $\mathcal{A}$ has a unique
kernel.

{\bf Claim (2).} Let $g:Y\rightarrow Z$ be a morphism in $\mathcal{H}$
such that $\underline{g}$ is an epimorphism. Then
$\underline{g}$ is a cokernel.

For such $g$ we form a triangle $X\s{f}{\rightarrow}
Y\s{g}{\rightarrow} Z\s{g'}{\rightarrow}X[1]$. We want to show that
$\underline{g}$ is the cokernel of $\underline{f}$.
 Let $h: Y\rightarrow M$ be a morphism with
 $\underline{h} \circ \underline{f}=0.$ Then we have some object $T_1\in
 \mathcal{T}$ and the following commutative diagram:

\[ \begin{CD}
X@>f>>Y@>g>>Z@>g'>>X[1]\\
@Vh_1 VV@VVh V @VVh_2V @VVh_3 V \\
T_1@>\delta >>M@>\delta _1>>Z'@>\delta _2>>T_1[1]
    \end{CD} \]

Since $\underline{g}$ is an epimorphism, by Theorem \ref{charmonoepi}, $g'$
factors through some object $T_2\in \mathcal{T}.$ It follows that
$h_3 \circ g'=0$, and then $h_2$ factors through $\delta _1$, i.e.
$h_2=\delta_1 \circ \delta_3$ for some morphism $\delta_3: Z\rightarrow
M$. As above it follows that $\underline{h}=\underline{\delta
_3} \circ \underline{g}$. Also the uniqueness of $\underline{\delta _3}$
is obtained in the same way as above. The corresponding
statement for monomorphisms can be shown dually.

This finishes the proof.
\end{proof}

\begin{cor} Let $\mathcal{H}$ and $\mathcal{T }$ be
as in Theorem \ref{charmonoepi}. and $M\s{f}{\rightarrow}
N\s{g}{\rightarrow}L\s{h}{\rightarrow}M[1]$ a triangle in
$\mathcal{H}$.

If $\underline{h}=0$ then $M\s{\underline{f}}{\rightarrow}
N\s{\underline{g}}{\rightarrow}L\rightarrow 0$ is exact in
$\mathcal{A}$.

If $\underline{h[-1]}=0$ then $0\rightarrow
M\s{\underline{f}}{\rightarrow} N\s{\underline{g}}{\rightarrow}L$ is
exact in $\mathcal{A}$.

Furthermore if
$\underline{h}=0=\underline{h[-1]},$ then $0\rightarrow
M\s{\underline{f}}{\rightarrow}
N\s{\underline{g}}{\rightarrow}L\rightarrow 0$ is an exact sequence
in $\mathcal{A}$.
\end{cor}

\section{More on the abelian quotient category}

In this section we will show that the quotient category modulo a tilting
subcategory is indeed the module
category of a certain endomorphism ring, which under mild additional
assumptions turns out to have various
strong properties, including being a Gorenstein algebra of Gorenstein
dimension at most one.
Several results in this section, especially in the first subsection,
generalize results of \cite{KR}, obtained there
under stronger assumptions. Results in the third subsection complement
results of \cite{BMR,BMRRT}.

\subsection{Endomorphisms algebras and Gorenstein property}

Here, and in the following we write $Ext(X,Y)$ for $Hom(X,Y[1])$.

\begin{lem}
Let $\mathcal{H}$ be a triangulated category,and $\mathcal{T }$ a full
subcategory of
$\mathcal{H}$ with $\Ext(\mathcal{T},\mathcal{T})=0$.
Then  $\mathcal{T}\cap \mathcal{T }[1]=\{0\}$.
\end{lem}

\begin{proof} This follows directly from the assumptions.
\end{proof}

\begin{prop} Let  $\mathcal{T }$ be a tilting subcategory of a triangulated
category  $\mathcal{H}$ and $\mathcal{A }$ the abelian quotient of
 $\mathcal{H}$ by $\mathcal{T}$. Then an object $M$ of $\mathcal{A}$ is a
projective object  if and only if $M\in \mathcal{T}[-1].$ Dually
 an object $N$ of $\mathcal{A}$  is an injective object  if and only if
$N\in\mathcal{T}[1].$
\end{prop}

\begin{proof} We prove the first statement only, the second one being dual.

Firstly we show that for any $T\in \mathcal{T}$ the shifted object $T[-1]$
is projective in $\mathcal{A}.$ For any epimorphism
 $X\s{\underline{f}}{\rightarrow }Y $
 in $\mathcal{A}$ and any morphism $\underline{g}: T[-1]\rightarrow Y$, let
$Z[-1]\rightarrow X\s{f}{\rightarrow} Y\s{h}{\rightarrow} Z$ be the triangle
into which $f$ is embedded. Since $\underline{f}$ is an epimorphism in
$\mathcal{A}$, the map $h$ factors through an object $T'$ of $\mathcal{T}$
by Theorem \ref{charmonoepi}. It follows that $h \circ g=0$ since
$Hom(T[-1],T')=0$. Then $g$ factors through $f$, hence $\underline{g}$
factors through $\underline{f}$.
This proves that $T[-1]$ is projective in $\mathcal{A}$.

Conversely assume $M$ is a projective object in $\mathcal{A}$.
By Lemma \ref{1-orthconditions}, there is a triangle
$T_1[-1]\rightarrow T_0[-1]\s{\sigma}{\rightarrow} M\rightarrow T_1$ with
$\sigma$ being a right $\mathcal{T}[-1]-$approximation of $M$. This yields
an exact sequence  $T_1[-1]\rightarrow T_0[-1]
\s{\underline{\sigma}}{\rightarrow} M\rightarrow 0$ with
$\underline{\sigma}$ being an epimorphism in $\mathcal{A}$.
So the sequence splits, hence $M\in \mathcal{T}[-1]$.
\end{proof}

The main result in this subsection is the following theorem, generalizing and
reproving in a different way a result in \cite{KR}.

By a category having enough projectives we mean that every object has
a projective cover.

An abelian category with enough projectives and enough injectives is
called {\em Gorenstein} if the full subcategory of projective
objects is covariantly finite and the full subcategory of inject
objects is contravariantly finite and there is an integer $d$ such
that all projectives are of injective dimension at most $d$ and all
injectives are of projective dimension at most $d$. The maximum of
the injective dimensions of projectives and the projective
dimensions of injectives is called {\em Gorenstein dimension} of the
category.

\begin{thm} Let $\mathcal{H}$ be a triangulated category, let
$\mathcal{T}$ be a tilting subcategory of $\mathcal{H}$ and let
$\mathcal{A}$ be the abelian quotient category of $\mathcal{H}$ by
$\mathcal{T}$. Then:
\label{Gorenstein}
\begin{enumerate}
\item The category $\mathcal{A}$ has enough projective objects.

\item The category $\mathcal{A}$ has enough injective objects.

\item The category $\mathcal{A}$ is Gorenstein of Gorenstein dimension at
most one.
\end{enumerate}
\end{thm}

\begin{proof}
We start by proving that any object $X$ of
$\mathcal{A}$ has a projective cover.
Now let $X\in \mathcal{H}$. There is a $\mathcal{T}[-1]-$approximation of $X$:
$T_1[-1]\s{f}{\rightarrow} X$, which is a morphism in the
triangle $T_2[-1]\rightarrow
T_1[-1]\s{f}{\rightarrow} X \s{g}{\rightarrow}T_2$ with
$T_2\in \mathcal{T }$. Thus we get a projective cover of $X$ in
$\mathcal{A}$: $\underline{f}: T[-1]\rightarrow X\rightarrow 0$ and a projective presentation of $X$:
  $T_2[-1]\rightarrow
T_1[-1]\s{f}{\rightarrow} X \rightarrow 0.$

Dually, injective objects in $\mathcal{A}$ are of the form $T[1]$ with $T\in
\mathcal{T}$, and any object has an injective envelope.

Furthermore, for any injective object $\T[1]$ in $\mathcal{A}$, we have a
 $\mathcal{T}[-1]-$approximation of $T[1]$:
$T_1[-1]\s{f}{\rightarrow}  T[1]$ since $\mathcal{T}[-1]$ is
contravariantly finite in $\mathcal{A}$. As before we have the
triangle $T\rightarrow T_2[-1]\rightarrow
T_1[-1]\s{f}{\rightarrow}T[1]\s{g}{\rightarrow}T_2$ with $T_2\in
\mathcal{T }$. By Theorem \ref{charmonoepi} there is an exact sequence:
$0\rightarrow T_2[-1]\rightarrow
T_1[-1]\s{\underline{f}}{\rightarrow}T[1]\s{\underline{g}}{\rightarrow}0$
which is a projective resolution of the
injective object $T[1]$ in $\mathcal{A}$.
Therefore proj.dim. $T[1]\leq 1.$  For a projective object $T[-1]$ of
$\mathcal{A}$, we have a triangle in $\mathcal{H}$:  $T\rightarrow
T[-1]\s{f}{\rightarrow}T_1[1]{\rightarrow}T_2[1]\rightarrow T$ with
$T_2\in \mathcal{T }$ and $f$ being a $\mathcal{T}[1]-$approximation
of $T[1]$. It follows that $0\rightarrow
T[-1]\s{\underline{f}}{\rightarrow}T_1[1]{\rightarrow}T_2[1]\rightarrow
0$ is an exact sequence in $\mathcal{A}$, which is an injective
resolution of the projective object $T[-1]$. Thus inj.dim$T[-1]\leq 1$.

Therefore $\mathcal{A}$ is an abelian category, which is Gorenstein of
Gorenstein dimension at most one.
\end{proof}

We denote by Mod$\mathcal{T}$ the category of modules over $\mathcal{T}$,
and by mod$\mathcal{T}$ the
subcategory of Mod$\mathcal{T}$ consisting of finitely presented modules.

As in \cite{KR} we get the following:

\begin{cor} Let $\mathcal{H}$ be a triangulated category
and $\mathcal{T }$ a tilting subcategory of $\mathcal{H}$. Then
$\mathcal{A}$ is equivalent to $\mod(\mathcal{T}[-1])$ as
abelian categories. \label{modulecat}
\end{cor}

\begin{proof} By Theorem \ref{Gorenstein}, $\mathcal{T}[-1]$ is
a full subcategory of the abelian category $\mathcal{A}$ consisting of
projective objects and
$\mathcal{A}$ has enough projectives and injectives. Therefore
$\mathcal{A}$ is equivalent to mod$(\mathcal{T}[-1])$, and the equivalence
preserves the exact structure.
\end{proof}

\begin{cor} Let $\mathcal{H}$ be a triangulated category
and $\mathcal{T }=add(T)$ a tilting subcategory of $\mathcal{H}$. Let
$A=End(T)$ be the endomorphism ring of $T$. Then $A$ is
Gorenstein of Gorenstein dimension at most one.
\end{cor}

Here, $A$ above may be an algebra of infinite dimension over a field $k$,
or it may just be a ring. If it is an artin algebra, then either it is
hereditary or its global dimension is infinite.
\medskip

 An abelian category with enough projectives and enough injectives is
called a {\em Frobenius category} if projective and injective objects
coincide.

\begin{prop} \label{Frobenius}
Let $\mathcal{H}$ be a triangulated category
and $\mathcal{T}$ a tilting subcategory. Then $\mathcal{A}$ is a Frobenius
category if and only if $\mathcal{T}=\mathcal{T}[2]$.
\end{prop}

\begin{proof} By the proof of Theorem
\ref{Gorenstein}, $\mathcal{A}$ is a Gorenstein
abelian category of Gorenstein dimension $\leq 1$ whose
projective objects are of the form $\mathcal{T}[-1]$ and
whose injective objects are of the form $\mathcal{T}[1]$. Then $\mathcal{A}$
is a Frobenius category if and only if $\mathcal{T}[1]=\mathcal{T}[-1]$
if and only if $\mathcal{T}=\mathcal{T}[2].$
\end{proof}

\subsection{Triangulated categories with Serre duality}

From now on, we assume that $\mathcal{H}$ is a $k-$linear triangulated
category with split idempotents and all Hom-spaces
of $\mathcal{H}$ are finite dimensional.  We also assume that $\mathcal{H}$
has Serre functor $\Sigma$ such that for all $X,Y$ there is a natural
isomorphism $Hom(X,Y) \simeq Hom(Y, \Sigma X)^{\ast}$, where $\ast$
denotes $k$-duality. Then $\mathcal{H}$ has
Auslander-Reiten triangles and $\Sigma=\tau [1],$ where $\tau$ is the
Auslander-Reiten translate. Without loss of generality, we may assume that
the AR-quiver of $\mathcal{H}$ has no loops. Indeed, triangulated
categories with loops in their AR quivers have been classified in
\cite{XZ}, Theorem 2.2.1. It turns out that in these categories
$\tau=[1]=id_{\mathcal{H}}$ and, obviously, $\mathcal{H}$ then has no
tilting subcategory.

\begin{prop}
Let $\mathcal{T}$ be a tilting subcategory of triangulated category
$\mathcal{H}$ and $\mathcal{A}$ the abelian quotient category of
$\mathcal{H}$ by $\mathcal{T}$. Then: \label{AR-properties}

\begin{enumerate}
\item $\mathcal{A}$ has source maps and sink maps.
In particular, the category $\mathcal{A}$ has AR-sequences.

\item $\mathcal{T}\cap \tau \mathcal{T }=\{ 0 \}$.

\item  There is equality $\tau^{-1}\mathcal{T}= \mathcal{T}[-1]$ i.e, $F\mathcal{T}=\mathcal{T}$, where $F=\tau ^{-1}[1].$

\end{enumerate}
\end{prop}

\begin{proof} It is routine to prove that the residue class of any sink (or
source) map in $\mathcal{H}$ is again a sink (or source, respectively)
 map in $\mathcal{A}$. Then $\mathcal{A}$ has sink maps and source maps,
 and it has AR-sequences.

Now we will prove the equality $\tau^{-1}\mathcal{T}= \mathcal{T}[-1]$.
For any projective object
$\tau^{-1}T$ with $T\in\mathcal{T}$, we have a
$\mathcal{T}[-1]-$approximation of $\tau^{-1}T$:
$T_1[-1]\s{f}{\rightarrow}\tau^{-1} T$ since $\mathcal{T}$, hence also
$\mathcal{T}[-1]$, are contravariantly finite in $\mathcal{A}$. As
before we have the triangle $T_2[-1]\rightarrow
T_1[-1]\s{f}{\rightarrow}\tau^{-1}T\s{g}{\rightarrow}T_2$ with
$T_2\in \mathcal{T }$. Since ${\Hom}(\tau^{-1}T, T_2)\cong
D\Ext^1(T_2,T)=0$, we have $g=0$. It follows that the triangle above
splits, i.e., $T_1[-1]\cong\tau^{-1}T\oplus \T_2[-1]$. This proves
$\tau^{-1}\mathcal{T}\subseteq \mathcal{T}[-1].$ A similar approximation
argument shows that $\tau^{-1}\mathcal{T}\supseteq \mathcal{T}[-1].$
Therefore $\tau^{-1}\mathcal{T}= \mathcal{T}[-1].$

For the proof of (2), we take $T\in \mathcal{T}\cap \tau \mathcal{T}$.
Then there exists
$T' \in \mathcal{T}$ such that  $T=\tau T'$.
Hence $\Hom (T,T)=\mbox{Hom}(T,\tau T')\cong D\mbox{Hom}(T',T[1])=0$,
and therefore $T'=T=0$.
\end{proof}

In the special case of cluster categories,
endomorphism rings of tilting objects
have been studied in \cite{BMR}. Assuming finite representation type,
a bijection has been shown to exist between the indecomposable
representations of the hereditary algebra and of the cluster tilted algebra.

\begin{prop} Let $\mathcal{H}$ be a $k-$linear
triangulated category over an algebraically closed field $k$ and
$\mathcal{T}_i =add(T_i)$ two tilting subcategories of $\mathcal{H}$,
for $i=1,\ 2$. Let $A_i=End(T_i)$ the endomorphism algebras of $T_i$. Then
$A_1$ and $A_2$ have the same representation type.
\end{prop}

\begin{proof} By Corollary \ref{modulecat},
$A_i-$mod$\thickapprox \mathcal{H}/ \mbox{add}T_i$ for $i=1,\ 2$.
Therefore  $A _1-\mbox{mod}/\mbox{add}T_2\thickapprox \mathcal{H}/
\mbox{add}(T_1\cup T_2)\approx A _2-\mbox{mod}/\mbox{add}T_1.$
Hence, ind$A_1$ is a finite set if and only if
ind$\mathcal{H}$ is a finite set. Thus $A_1$ is of finite type if
and only if $A_2$ is so.
Moreover, by \cite{Kr}, $A _1-\mbox{mod}$ is  wild if and only if  $A
_2-\mbox{mod}$ is wild.  Therefore, by tame-wild dichotomy, $A_1$ and
$A_2$ have the same representation type.
\end{proof}

\begin{thm} \label{image-1-orthogonal}
Let $\mathcal{T}$ be a tilting subcategory of
$\mathcal{H}$  and $\mathcal{A}$ the abelian quotient. If
$\mathcal{C}$ is a $1-$orthogonal subcategory of $\mathcal{H}$, i.e.
$Ext^1_{\mathcal{H}}(\mathcal{C},\mathcal{C})=0$,
then its image in $\mathcal{A}$ is a $1-$orthogonal subcategory of
$\mathcal{A}$, i.e. $Ext^1_{\mathcal{A}}(\mathcal{C},\mathcal{C})=0$.
\end{thm}

\begin{proof} Let $X, \ X_1\in \mathcal{C}$ such that $X_1$ has no
direct summands in $\mathcal{T}$. We will prove that
$\Ext^1_{\mathcal{A}}(X_1, X)=0$, i.e., any short exact sequence
$0\rightarrow
X\s{\underline{f}}{\rightarrow}M\s{\underline{g}}{\rightarrow}X_1\rightarrow
0$ in $\mathcal{A}$ splits.  Lifting the morphism
$\underline{f}: X\rightarrow M$ in $\mathcal{A}$ to a morphism $f$ in
$\mathcal{H}$, we get a triangle $N[-1]\s{f_1}{\rightarrow}
X\s{f}{\rightarrow}M\s{f_2}{\rightarrow}N\s{f_3}\rightarrow X[1]$.
Since $\underline{f_1} =0$,  $\underline{f}$ is a monomorphism.

From our construction of the cokernel of a monomorphism in the proof of
\ref{charmonoepi}, we get the
following commutative diagram:

$$ \begin{array}{lclclclcl}
&&&&&&T_1Y&=&T_1\\
&&&&&&\downarrow  &&\downarrow \\
X'[-1]&\s{g_3}{\rightarrow} &T_0[-1]&\s{g_2}{\rightarrow}&M&\s{g}{\rightarrow}& X'&\s{g_1}{\rightarrow}&T_0\\
\downarrow h_4&&\downarrow h_3&&\parallel &&\downarrow h_2&&\downarrow h_1\\
N[-1]&\s{f_1}{\rightarrow}&X&\s{f}{\rightarrow}&M&\s{f_2}{\rightarrow}&N&\s{f_3}{\rightarrow}&X[1]
\end{array}$$

Since $X'\cong X_1$ in $\mathcal{A}$ we have that $X'\cong X_1\oplus
T'$ for some $T'\in \mathcal{T}.$ Then $h_4$ can be written as
$h_4=(h_5, h_6): X_1[-1]\oplus T'[-1]\rightarrow N[-1]$. It follows that
$f_1 \circ h_4=(f_1 \circ h_5, f_1 \circ h_6)$ where
$f_1 \circ h_5\in \mbox{Hom}(X_1[-1], X)=0$
and $f_1 \circ h_6=0$, the latter because $f_1$ factors through an
object in $\mathcal{T}$ and thus $f_1 \circ h_6$ factors through a map
from some $T''[-1]$ to some $T'''$, which by assumption is zero.
Therefore $h_3 \circ g_3=0$ and so there exists
 a morphism $\sigma :M\rightarrow X$ such that $h_3=\sigma \circ g_2$. It
 follows that $(1-f\sigma) \circ g_2=0.$ Hence
there exists a morphism $\rho
 :X'\rightarrow M$ such that $1-f \circ \sigma=\rho \circ g$ i.e.
$1=f \circ \sigma+\rho \circ g$.
By passing this equality to the quotient category
 $\mathcal{A}$, we get that
 $1=\underline{f} \circ \underline{\sigma}+\underline{\rho} \circ
 \underline{g}.$  Here, $e_1=\underline{f}\underline{\sigma}$ and
$e_2=\underline{\rho}
 \underline{g}$ are orthogonal idempotents of End$_{\mathcal{A}}$.
 Then $M\cong e_1M\oplus e_2M$, and $e_1M=\underline{f}
\underline{\sigma}M\cong \underline{\sigma}M$ and
$e_2M=\underline{\rho} \underline{g}M\cong \underline{\rho}X_1$. So $e_1M$
is a subobject of $X$ and $e_2M$ is an image of $X_1$. Since  $0\rightarrow
X\s{\underline{f}}{\rightarrow}M\s{\underline{g}}{\rightarrow}X_1\rightarrow
0$ is an exact sequence, we obtain that $\underline{\sigma}M\cong X$ and
$\underline{\rho}X_1\cong X_1$ by computing their lengths.
 Therefore the exact sequence $0\rightarrow
X\s{\underline{f}}{\rightarrow}M\s{\underline{g}}{\rightarrow}X_1\rightarrow
0$ in $\mathcal{A}$ splits. This finishes the proof.
\end{proof}
\medskip

\begin{cor} Under the same assumptions as in Theorem \ref{image-1-orthogonal},
an indecomposable $1-$orthogonal object $C$ in $\mathcal{H}$ (that is,
$Ext^1_{\mathcal{H}}(C,C)=0$), which does not belong to
$\mathcal{T}$ is an $1-$orthogonal indecomposable object in $\mathcal{A}$.
\end{cor}

Such $1-$orthogonal objects sometimes are also called exceptional
objects.

\begin{prop}
Let $\mathcal{H}$ be a triangulated category
and $\mathcal{T}$ a tilting subcategory. Then $mod\mathcal{T}[-1]$ is a
Frobenius category if and only if $\Sigma\mathcal{T}=\mathcal{T}$.
\end{prop}

\begin{proof} By Proposition \ref{AR-properties}, $\tau^{-1}\mathcal{T}
=\mathcal{T}[-1]$. Then by Proposition \ref{Frobenius},
$mod\mathcal{T}[-1]$  is a Frobenius category if and only if
$\mathcal{T}=\mathcal{T}[2]$ if and only if  $\Sigma\mathcal{T}=\mathcal{T}$.
\end{proof}

\subsection{Cluster categories}

Cluster categories are the motivating example for our results. These
categories have been introduced in \cite{BMRRT}, and in \cite{CCS1} in
the case of type $A_n$, in order to connect the cluster algebras defined
by Fomin and Zelevinsky \cite{FZ1,FZ2} (see also the survey  \cite{FZ3} on
cluster algebras), with representation theory of
algebras. The cluster variables of Fomin and Zelevinsky correspond to
indecomposable exceptional objects in cluster categories, and clusters
correspond to tilting objects, that is, to maximal $1-$orthogonal
subcategories \cite{I1}, which play a crucial role also in our more general
framework. For recent developments on cluster tilting, we refer
to the survey papers  \cite{BM, Rin}.

Recall that cluster categories by definition are orbit categories
$D^{b}(H)/F$ of derived categories $D^{b}(H)$ (of a hereditary
category $H$) by an automorphism group generated by $F=\tau ^{-1}[1]$
where $\tau$ is the Auslander-Reiten
translate in $D^{b}(H)$, and $[1]$ is the shift functor of
$D^{b}(H)$. Cluster categories are triangulated categories by [K] and they
form examples of Calabi-Yau triangulated categories of CY-dimension $2$ as
studied in \cite{KR}.

Of particular interest are the endomorphism algebras of tilting objects.
These provide, or are expected to provide, essential information on cluster
variables and clusters, see \cite{CCS1, CCS2}. Moreover,
by \cite{BMR,BMRRT,KR}
quotients of cluster categories or Calabi-Yau categories of CY-dimension
two modulo tilting objects are equivalent to the module category of the
corresponding endomorphism algebra. Our main theorem \ref{mainthm} puts these
results into a more general context. Moreover, several results we prove in
the present section are direct generalisations of results on cluster or
Calabi-Yau categories \cite{BMR,KR,Z2}, for instance on representation types
or on the Gorenstein property.

We add another result in this special situation:

\begin{cor} Let $T$ be a tilting object of a cluster
category $\mathcal{C(H)}$ and $A$ the cluster tilted algebra. If
$A$ is hereditary and $T'$ is a tilting object in $\mathcal{C(H)}$ with
$\add(T)\cap\mbox{add}(T')=\{ 0\},$ then
 $T'$ is a tilting module in $A-$mod.
\end{cor}

\begin{proof} By Theorem \ref{image-1-orthogonal}, $T'$ is a partial tilting
$A-$module. It is a tilting module since the number of non-isomorphic
indecomposable summands of $T'$ and of $T[-1]$ is the correct one.
\end{proof}
\medskip

In the following, we apply our results to cluster categories. Let $H$ be a
hereditary algebra and $F=\tau^{-1}[1]$. $F$ is an automorphism of
$D^b(H)$. The cluster category  $\mathcal{C}(H)= D^b(\mathcal{H})/F$ is
 a triangulated category, and the projection  $\pi: D^b(H)\longrightarrow
\mathcal{C}(H)$ is a triangle functor. Now we show that it induces a
covering functor of cluster tilted algebras.

\begin{thm} Let $\mathcal{T}$ be a tilting subcategory of $D^b(H)$ and
$\pi: D^b(H)\longrightarrow \mathcal{C}(H)$ the projection.
 Then:
\begin{enumerate}
\item The restriction of $\pi$ to $\mathcal{T}[-1]$ is a Galois covering of
the cluster tilted algebra $\pi(\mathcal{T}[-1])$.
\item  The projection $\pi$ induces a covering functor from
$mod(\mathcal{T}[-1])$ to $mod(\pi(\mathcal{T}[-1]))$.
\end{enumerate}
\end{thm}

This is closely related to results in section two of \cite{BMRRT}, where
Ext-configurations and tilting sets are studied. In particular,
Propositions 2.1 and 2.2. there compare properties relevant to tilting
in the derived category and in the cluster category.

Before we prove the Theorem, we first show a lemma:

\begin{lem} $\mathcal{T}$ is a tilting subcategory of  $D^b(H)$ if and only
if $\mathcal{T}=\pi ^{-1}(\pi(\mathcal{T}))$ and
$\pi(\mathcal{T})$ is a tilting subcategory of $\mathcal{C}(H)$.
\label{tiltinginverse}
\end{lem}

\begin{proof} For a subcategory $\mathcal{T}$ of $D^b(H)$ with
$\mathcal{T}=\pi ^{-1}(\pi(\mathcal{T}))$, it is easy to prove that
 $\mathcal{T}$ is contravariantly finite in $D^b(H)$ if and only if
$\pi(\mathcal{T})$ is so in $\mathcal{C}(H).$

Suppose $\mathcal{T}$ is a tilting subcategory of $D^b(H)$. Then
$F\mathcal{T}=\mathcal{T}$ by Proposition \ref{AR-properties}. We
denote by $\mathcal{T}'$ the intersection of $\mathcal{T}$ with the
additive subcategory $\mathcal{C}'$  generated by all $H-$modules
as stalk compleses of degree $0$ together with
$H[1]$. Then we have that $\mathcal{T}=\{ F^n(\mathcal{T}')| n\in
\mathbf{Z} \}.$ Now $\pi(\mathcal{T})=\pi(\mathcal{T}')$, denoted by
$\mathcal{T}_1$. For any pair of objects
$\tilde{T}_1, \tilde{T}_2\in \mathcal{T}_1$, there are $T_1, T_2\in
\mathcal{T}'$ such that $\tilde{T}_1 =\pi (T_1), \tilde{T}_2=\pi (T_2)$.
 Then  $Ext^1(\tilde{T}_1, \tilde{T}_2)=Hom(\tilde{T}_1, \tilde{T}_2[1])
\cong \oplus _{n\in Z}Hom_{D^b(H)}(T_1, F^nT_2[1])
= \oplus _{n\in Z}Ext^1_{D^b(H)}(T_1, F^nT_2)=0.$ If there are
indecomposable objects $\tilde{X}=\pi (X)\in \mathcal{C}(H)$ with $X\in
\mathcal{C}'$ satisfying
$Ext^1(\mathcal{T}_1,\tilde{X})=0$, then $Ext^1(F^n\mathcal{T}', X)=0$ for
any $n$, and then $Ext^1(\mathcal{T}, X)=0$. Hence $X\in \mathcal{T}$ by
$\mathcal{T}$ being a tilting subcategory. Thus $\tilde{X}\in
\mathcal{T}_1.$ This proves that the image $\mathcal{T}_1$  of $\mathcal{T}$
under $\pi$ is a tilting subcategory of $\mathcal{C}(H).$

Conversely, from $\mathcal{T}=\pi^{-1}(\mathcal{T}_1)$, we get
$F(\mathcal{T})= \mathcal{T}$. As above we denote
by $\mathcal{T}'$ the intersection of $\mathcal{T}$ with the additive
subcategory $\mathcal{C}'$  generated by all $H-$modules
as stalk compleses of degree $0$ together with
$H[1]$.  Since $Ext^1(\mathcal{T}_1, \mathcal{T}_1)
\cong\oplus _{n\in Z} Ext^1(\mathcal{T}',F^n\mathcal{T}')=0$, we have that
$Ext^1(F^m\mathcal{T}',F^n\mathcal{T}')\cong Ext^1(\mathcal{T}', F^{n-m}
\mathcal{T}')=0.$
This proves that $\mathcal{T}$ is an orthogonal subcategory. Now if $X\in
D^b(H)$ satisfies $Ext^1(\mathcal{T}, X)=0$,
then $Ext^1(\mathcal{T}_1, \tilde{X})=0$. It follows that $\tilde{X}\in
\mathcal{T}_1$, hence $X\in \mathcal{T}$.
Similarly, if $X\in D^b(H)$ satisfies $Ext^1(X,\mathcal{T})=0$, then
$X\in \mathcal{T}$.
\end{proof}

Now we are ready to give the proof of the theorem.

\begin{proof}
(1). By Lemma \ref{tiltinginverse}, $\pi (\mathcal{T} )$ is a
cluster tilting object in $\mathcal{C}(H).$ The projection $\pi $ sends
$\mathcal{T}[-1]$ to $\pi(\mathcal{T}[-1])$, which is equivalent to
the cluster tilted algebra $\pi(\mathcal{T})[-1]$ since $\pi$ is a triangle functor.
Thus $\pi_{\mathcal{T}[-1]}: \mathcal{T}[-1]\longrightarrow
\pi(\mathcal{T})[-1]$ is a Galois covering with Galois group generated
by $F$.

(2). By Theorem \ref{Gorenstein} and Corollary \ref{modulecat} (or by
\cite{BMR}, \cite{Z1}, \cite{KR}) there are equivalences
$D^b(H)/\mathcal{T}\cong mod(\mathcal{T}[-1])$ and
$\mathcal{C}(H)/(\pi(\mathcal{T}))\cong mod(\pi(\mathcal{T})[-1])$.
We define the induced functor $\bar{\pi}$ as follows:
$\bar{\pi}(X):= \pi (X)$  for any object $X\in D^b(H)/\mathcal{T}, $ and
$\bar{\pi}(\underline{f}):=\underline{\pi (f)}$ for any morphism
$\underline{f}: X\rightarrow Y$ in $D^b(H)/\mathcal{T}$. Clearly
$\bar{\pi}$ is well-defined and makes the following diagram
commutative:

\[ \begin{CD}
D^b(H)@>\pi>>\mathcal{C}(H)\\
@Vq_1 VV@VVq_2 V \\
D^b(H)/\mathcal{T}@>\bar{\pi} >>\mathcal{C}(H)/\pi(\mathcal{T}).
    \end{CD} \]

Then $\bar{\pi}$ is a covering functor from $D^b(H)/\mathcal{T}$ to
$\mathcal{C}/\pi(\mathcal{T})$, i.e, it is a covering functor
from $mod(\mathcal{T}[-1])$ to $mod(\pi(\mathcal{T}))$.
\end{proof}

\subsection{Self-injective algebras}

Stable module categories of self-injective algebras are triangulated
categories with Serre functor. Preprojective algebras
and group algebras of finite groups are examples of self-injective algebras.

\begin{prop} Let $A$ be a self-injective finite dimensional  algebra and
  $M$ an $A-$module.
 Then $M$ is a maximal $1-$orthogonal module if and only if add$(M)$ is a
tilting subcategory of $A-\underline{\mbox{mod}}$.
\end{prop}

\begin{proof} We note that  $Ext^1_A(X,Y)\cong
\underline{\mbox{Hom}}(X,\Omega ^{-1} Y)$, for any $A-$modules $X, Y .$ Then
 $Ext^1_A(X,Y)\cong Ext^1_{A-\underline{mod}}(X, Y).$ It follows that $M$ is
a maximal $1-$orthogonal module if and only if add$(M)$
 is a tilting subcategory of $A-\underline{\mbox{mod}}$.
\end{proof}

A maximal 1-orthogonal module over a self-injective algebra must contain a
projective generator. Hence we get:

\begin{cor}  Let $A$ be a self-injective finite dimensional  algebra and
$M$ a maximal $1-$orthogonal module. Then
$A-\mbox{mod}/\mbox{add}M$ is again an abelian category.

\end{cor}

\subsection{Other examples}

The following examples indicate that our examples cover not only cluster
categories, but also some stable categories. Moreover, we also
cover situations not of Calabi Yau dimension two.

\begin{enumerate}
\item Let $\mathcal{H} = A-\underline{mod}$ be the stable category of
the self-injective algebra $A=kQ/I$ given by the quiver $Q$:\\
\vspace*{-1cm}
\begin{center}
 \setlength{\unitlength}{0.61cm}
 \begin{picture}(5,4)
 \put(0,2){a}\put(0.4,2.2){$\circ$}
\put(3,2.2){$\circ$}\put(3.4,2){b} \put(0.8,2.5){\vector(3,0){2}}
\put(2.8,2.2){\vector(-1,0){2}}
 \put(1.7,2.7){$\alpha$}\put(1.7,1.5){$\beta$}
 \end{picture}
 \end{center}
\vspace*{-1cm}
 modulo the relations
 $\alpha \beta\alpha \beta, \  \beta\alpha \beta\alpha $.
 $\mathcal{H}$ is not of CY-dimension $2$.

The following is the Auslander Reiten quiver of $A-mod$ (the first
and the last column have to be identified). Deleting the top row gives
the Auslander Reiten quiver of $A-\underline{mod}$.

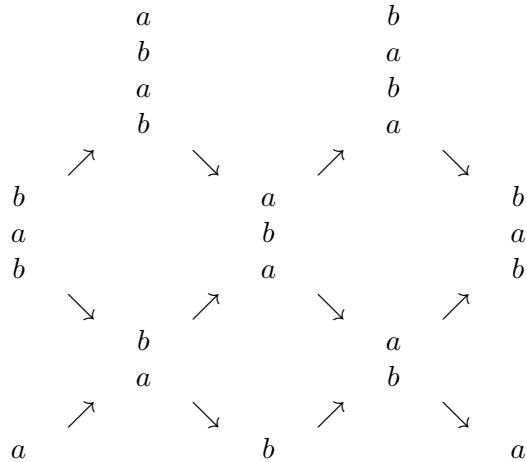
\newcommand{\Mabab}{\begin{array}{c} a \\ b \\ a \\ b \end{array}}
\newcommand{\Mbaba}{\begin{array}{c} b \\ a \\ b \\ a \end{array}}
\newcommand{\Maba}{\begin{array}{c} a \\ b \\ a \end{array}}
\newcommand{\Mbab}{\begin{array}{c} b \\ a \\ b \end{array}}
\newcommand{\Mab}{\begin{array}{c} a \\ b \end{array}}
\newcommand{\Mba}{\begin{array}{c} b \\ a \end{array}}

$$\begin{array}{ccccccccc}
&& \Mabab &&&& \Mbaba && \\
& \nearrow && \searrow && \nearrow && \searrow & \\
\Mbab &&&& \Maba &&&& \Mbab \\
& \searrow && \nearrow && \searrow && \nearrow & \\
&& \Mba &&&& \Mab && \\
& \nearrow && \searrow && \nearrow && \searrow & \\
a &&&& b &&&& a
\end{array}$$

We choose $\mathcal{T}$
to be the subcategory $\add(M)$ where $M$ is the direct sum of
the simple module $L(a)$ and the indecomposable module
$\Maba$ with top and socle isomorphic to $L(a)$ and with length $3$.
Then $\mathcal{T}$ is  a tilting subcategory of $\mathcal{H}$.
The quotient category of
$\mathcal{H}$ by this tilting subcategory $\mathcal{T}$
is equivalent to the module
category of the endomorphism algebra $B$ of $M$. Here $B$ is given
by the same quiver with relation $\alpha \beta, \  \beta\alpha$.

The Auslander Reiten quiver of the quotient category is as follows:

$$\begin{array}{ccccccccc}
&& \Mba &&&& \Mab && \\
& \nearrow && \searrow && \nearrow && \searrow & \\
\Mbab &&&& b &&&& \Mbab
\end{array}$$

Again the first and the last column are identified.

Note that north-east arrows denote epimorphisms, while south-east
arrows denote monomorphisms.
\medskip

\item Let $\mathcal{H}$ be the (bounded) derived category of
hereditary algebra $A$, where $A$ is the path algebra of the quiver: \\
\vspace*{-1cm}
\begin{center}
 \setlength{\unitlength}{0.61cm}
 \begin{picture}(5,4)
 \put(0,2){a}\put(0.4,2.2){$\circ$}
\put(3,2.2){$\circ$}\put(3.4,2){b}
\put(6,2.2){$\circ$}\put(6.4,2){c}
\put(0.8,2.5){\vector(3,0){2}}
\put(3.8,2.5){\vector(3,0){2}}
 \end{picture}
 \end{center}
\vspace*{-1cm}

Let $P_a, P_b, P_c$ be the indecomposable projective modules of $A$
with simple top $L(a), \ L(b), \ L(c)$, respectively.
If we take $\mathcal{T}$ to
be the subcategory generated by $\{ \tau^{-n}P_a[n],
\tau^{-n}P_b[n], \tau^{-n}P_c[n]\mid n \in Z\},$ then $\mathcal{T}$
is a tilting subcategory of $\mathcal{H}$ and
$D^b(A)/\mathcal{T}\cong \oplus _{i\in Z}A_i$ where $A_i\cong A$ for
any $i$.

 If we take $\mathcal{T}'$ to
be the subcategory generated by $\{ \tau^{-n}P_a[n],
\tau^{-n}L(c)[n], \tau^{-n}P_c[n]\mid n \in Z\},$ then
$\mathcal{T}'$ is also a tilting subcategory of $\mathcal{H}$ and
$D^b(A)/\mathcal{T}\cong B$ where $B$ is the locally finite path
algebra of the quiver
$$A^{\infty}_{\infty}:   \cdots \circ \longrightarrow \circ
\longrightarrow \circ \longrightarrow\cdots $$ with
$\underline{r}^2=0.$
\end{enumerate}

\section{Partial converse}

In this section, we discuss potential converse results to Theorem
\ref{mainthm}. Obviously, the direct converse does not hold true. A
trivial counterexample comes from the trivial category - with a zero
object only - being abelian. A more interesting counterexample is given
at the end of this section; a non-trivial abelian quotient category
obtained by factoring out a subcategory that is not tilting.

Other counterexamples can be obtained by starting with the derived
module category of a finite dimensional path algebra of a quiver
and factoring
out all preprojective and preinjective components. For example, start
with the tame Kronecker algebra $H$ (over an algebraically closed field
$k$), which is derived equivalent to the category of coherent sheaves over
the projective line. The indecomposable objects in
${\mathcal H} = D^b(H-mod)$ are shifts
of indecomposable modules, which are either preprojective or regular or
preinjective. Let $T$ be the full subcategory generated by all sums of
shifts of preprojective or preinjective modules (that is, of torsionfree
sheaves). Then the quotient category
${\mathcal H}/T$ has as objects all shifts of
regular modules. In the quotient category there are no maps between regular
objects in different degrees, since the extensions between regulars existing
in the module category are maps factoring through injective objects.
And in each fixed degree, the category of
regular objects decomposes into blocks, called tubes,
each of which is equivalent
to the category of finite dimensional modules over a power series ring in
one variable. Therefore, the quotient ${\mathcal H}/T$ also decomposes
into blocks of this type. Such a tube is an abelian category without
projective or injective objects.

\begin{thm} Let $\mathcal{H}$ be a triangulated category
and $\mathcal{T }$ a full subcategory of $\mathcal{H}$. Suppose that
$\mathcal{A}$ is an abelian category (with induced structure).
Then the following conditions are equivalent:

\begin{enumerate}
\item ${\Hom}(T,T'[1])=0$ for any $T, T'\in \mathcal{T }$.

\item $\mathcal{T }\cap\tau \mathcal{T }=\{ 0 \}$ and for any triangle
$Z[1]\s{h}{\rightarrow}X\s{f}{\rightarrow}Y\s{g}{\rightarrow}Z,$ if
$\underline{h}=0$, then the map
$\underline{f}$ is a monomorphism  in $\mathcal{A}$ .

\item $\mathcal{T }\cap\tau \mathcal{T }=\{ 0 \}$  and for any triangle
$Z[1]\s{h}{\rightarrow}X\s{f}{\rightarrow}Y\s{g}{\rightarrow}Z$,  if
$\underline{g}=0$, then the map
$\underline{f}$ is an epimorphism in $\mathcal{A}$ .
\end{enumerate}
\end{thm}

Since the abelian structure is induced from the triangulated one,
Auslander Reiten triangles in the quotient become Auslander
Reiten sequences (if non-trivial).

\begin{proof} $(2)\Rightarrow (3)$ Suppose $s: Y\rightarrow M$
satisfies $\underline{s} \circ \underline{f} = 0$, hence
$s \circ f$ factors through $\mathcal{T}$: Then there exists
the following commutative diagram, with $T\in \mathcal{T}$:

\[ \begin{CD}
X@>f>> Y@>g>>Z@>\sigma>>X[1]\\
@Vs'VV@VVs V @VVs''V @VV V \\
T@>f'>>M@>g'>>Z'@>>>T[1]
    \end{CD} \]

Then $g'\circ s=s'' \circ g$, so $\underline{g'} \circ \underline{s}
=\underline{s''} \circ \underline{g}=0$. By (2),
$\underline{g'}$ is a monomorphism, since $\underline{f'} = 0$.
This implies that $\underline{s}=0$.

$(3) \Rightarrow (2)$ Similar to the argument for the converse implication.

$(1) \Rightarrow (2)$  The statement is part of Theorem
\ref{charmonoepi}.

$(2)\Rightarrow (1)$  Suppose there are indecomposable objects $T,
T'\in \mathcal{T } $ with ${\Hom}(T,T'[1])\not=0.$ Then there is a
non-zero morphism $\sigma : T'\rightarrow \tau T$ since
$0\not=\mbox{Hom}(T,T'[1])\cong D\mbox{Hom}(T'[1],\tau T[1])\cong
D\mbox{Hom}(T',\tau T).$ Let
 $T[-1]\s{h}{\rightarrow}\tau
 T\s{g}{\rightarrow}M\s{}{\rightarrow}T$ be the AR-triangle ending at $T$
and $T'\s{\sigma}{\rightarrow}\tau
T\s{f}{\rightarrow}N\s{}{\rightarrow}T'[1]$
the triangle into which $\sigma$ is
 embedded. Then $\underline{g}$ is an epimorphism
 and $\underline{f}$ is a monomorphism by the condition (2). Moreover,
there is a morphism $h: M\rightarrow N$ such that $f=h \circ g$ since
  $f$ is not split. Then $\underline{g}$ is also a monomorphism, hence it is an isomorphism in $\mathcal{A}$.  This is a
  contradiction to $\underline{g}$ being non-zero and a source map.
\end{proof}

\begin{thm}
Let $\mathcal{H}$ be a triangulated category and $\mathcal{T }$ a
full subcategory of $\mathcal{H}$ with $\Ext^1(\mathcal{T
},\mathcal{T })=0$. Suppose that $\mathcal{A}$ is an abelian
category (with induced structure). Then for any $X\in \mathcal{H}$,
if $\Ext^1(X,\mathcal{T })=0$ and $\Ext^1(\mathcal{T },X)=0$, then
$X\in \mathcal{T}$.
\end{thm}

\begin{proof} Suppose $X$ is an indecomposable object satisfying
$ßExt^1(X,\mathcal{T })=0$ and $\Ext^1(\mathcal{T },X)=0$.
Assume $X\notin \mathcal{T }$.
Since ${\Hom}(X,\tau X[1])\cong D\mbox{Hom}(X,X)\not= 0$, $\tau X\notin
\mathcal{T }$. Let $\tau X\s{f}{\rightarrow} M\s{g}{\rightarrow}
X\s{h}{\rightarrow} \tau X[1]$ be the AR-triangle ending at $X.$

{\bf Claim:} $M\in \mathcal{T }$. \\ Otherwise, $0\rightarrow \tau
X\s{\underline{f}}{\rightarrow} M\s{\underline{g}}{\rightarrow}
X\rightarrow0$ is the AR-sequence ending at $X$ in the abelian
category $\mathcal{A}$.  Then $\underline{g}$ is an epimorphism in
$\mathcal{A}.$  Hence $h[-1]: X[-1]\rightarrow \tau X$ factors
through some object of $\mathcal{T }$, i.e., there are $T\in
\mathcal{T }$ and morphisms $h_1: X[-1]\rightarrow T$ and $h_2:
T\rightarrow \tau X$ such that $h[-1]=h_2 \circ h_1$. Hence $h[-1]=0$ since
$h_1\in \mbox{Hom}(X[-1], T)\cong \mbox{Hom}(X, T[1])=0.$ Thus
$h=0$, a contradiction.

Now let
$X\s{f_1}{\rightarrow}\tau^{-1} M\s{g_1}{\rightarrow} \tau^{-1}
X\s{h_1}{\rightarrow}  X[1]$ be the AR-triangle starting at $X.$

{\bf Claim:} $\tau ^{-1}X\notin \mathcal{T }$. Otherwise, for any
 $T\in \mathcal{T }$, we have that $0=\mbox{Hom}(\tau ^{-1}X,
 T[1])\cong\mbox{Hom}(X,\tau T[1])\cong D\mbox{Hom}(T, X)$. This is
 a contradiction to the AR-triangle ending at $X$ and having middle
term $M\in\mathcal{T}$.

Therefore we have the following
AR-sequence in $\mathcal{A}$ starting at $X$:
$0\rightarrow X\s{\underline{f_1}}{\rightarrow}\tau ^{-1}
M\s{\underline{g_1}}{\rightarrow} \tau ^{-1}X\rightarrow0$. In
particular, $\underline{g_1}$ is an epimorphism and thus
$\underline{h_1}=0$, i.e., there are
morphisms $\tau ^{-1}X\s{h_2}{\rightarrow}T\s{h_3}{\rightarrow}X[1]$
with $h_1=h_3 \circ h_2.$ Thus $h_1=0$ since $h_3\in \mbox{Hom}(T,X[1])=0$,
and it follows that $h=0$. Hence the AR triangle splits,
a contradiction. So we get $X\in \mathcal{T }$.
\end{proof}
\medskip

{\bf Example.} The following example gives a situation not covered by
our results.

Let $A=kQ/I$ be the self-injective algebra given by the quiver $Q$

\vspace*{-1cm}
\begin{center}
 \setlength{\unitlength}{0.61cm}
 \begin{picture}(5,4)
 \put(0,2){a}\put(0.4,2.2){$\circ$}
\put(3,2.2){$\circ$}\put(3.4,2){b} \put(0.8,2.5){\vector(3,0){2}}
\put(2.8,2.2){\vector(-1,0){2}}
 \put(1.7,2.7){$\alpha$}\put(1.7,1.5){$\beta$}
 \end{picture}
 \end{center}
\vspace*{-1cm}
and the relations $\alpha \beta\alpha, \  \beta\alpha \beta $.

\newcommand{\Maba}{\begin{array}{c} a \\ b \\ a \end{array}}
\newcommand{\Mbab}{\begin{array}{c} b \\ a \\ b \end{array}}
\newcommand{\Mab}{\begin{array}{c} a \\ b \end{array}}
\newcommand{\Mba}{\begin{array}{c} b \\ a \end{array}}

The Auslander Reiten quiver of $A-mod$ looks as follows:
$$\begin{array}{ccccccccc}
\Mbab &&&& \Maba &&&& \Mbab \\
& \searrow && \nearrow && \searrow && \nearrow & \\
&& \Mba &&&& \Mab && \\
& \nearrow && \searrow && \nearrow && \searrow & \\
a &&&& b &&&& a
\end{array}$$
Here, the first and the last column are identified.

Deleting the first row produces the Auslander Reiten quiver of
the stable category $A-\underline{mod}$:

$$\begin{array}{ccccccccc}
&& \Mba &&&& \Mab && \\
& \nearrow && \searrow && \nearrow && \searrow & \\
a &&&& b &&&& a
\end{array}$$

This stable category ${\mathcal H} = A-\underline{mod}$ of $A$
has no tilting subcategory. Indeed, including any of the four
indecomposable objects into $\mathcal{T}$ leaves us with the problem
that one of the maximality conditions forces us to include another
object, since it has no extensions with the first one, in one direction;
but then there are always extensions in the other direction, thus spoiling
another defining condition.

But ${\mathcal H} = A-\underline{mod}$ does have non-trivial abelian
quotient categories. For example, choosing
$\mathcal{T}$ to be the subcategory $add(L(a))$ of $\mathcal{H}$,
it is not difficult to check that $\mathcal{H}/\mathcal{T}$ is an
abelian category.

The Auslander Reiten quiver of $\mathcal{H}/\mathcal{T}$ is:

$$\begin{array}{ccccc}
\Mba &&&& \Mab  \\
& \searrow && \nearrow & \\
&& b &&
\end{array}$$

Note that the arrow pointing south-east represents a monomorphism, while
the arrow pointing north-east represents an epimorphism in the abelian
quotient category. There is a projective object, which is not
of the form $\mathcal{T}[-1]$.

\end{document}